\documentclass[10pt]{article}
\usepackage{amsmath,amssymb,amsfonts, euscript}

\newtheorem{theorem}{Theorem}
\newtheorem{lemma}[theorem]{Lemma}
\newtheorem{corollary}[theorem]{Corollary}
\newtheorem{proposition}[theorem]{Proposition}

\newtheorem{remark}[theorem]{Remark}
\newtheorem{algorithm}[theorem]{Algorithm}
\begin{document}

\title{Cohen-Macaulayness of bipartite graphs, revisited}
\author{Rashid Zaare-Nahandi\\ Institute for Advanced Studies in Basic Sciences,
Zanjan 45195, Iran \\ E-mail: rashidzn@iasbs.ac.ir}

\date{}

\maketitle

\begin{abstract}
Cohen-Macaulayness of bipartite graphs is investigated by several mathematicians and has been characterized combinatorially. In this note, we give some different combinatorial conditions for a bipartite graph which are equal to Cohen-Macaulayness of the graphs. Conditions in the previous works are depending on an appropriate ordering on vertices of the graph. The conditions presented in this paper are not depending to any ordering. Finally, we present a fast algorithm to check Cohen-Macaulayness of a given bipartite graph.


\noindent{\bf 2010 MR Subject Classification:} 13F55, 05C25, 05E45.
\end{abstract}

Characterization and classification of Cohen-Macaulay graphs, specially bipartite graphs, have been extensively studied in the last decades. For instance, see \cite{EsVil}, \cite{HH}, \cite{HHZ}, \cite{Vill1} and \cite{HYZ}.
Complete prerequisites for the subject are nicely written in the mentioned references and \cite{St}.  To make this note self-contained, we review here some basic definitions.

Through out this paper,  $G$ is a finite simple graph with no any vertex of degree zero. For two vertices $v$ and $w$ which are adjacent in $G$, we write $v\sim w$. The set of all vertices of $G$ adjacent to a vertex $v$ is denoted by $N(v)$.
A subset $P$ of the set of edges is called a perfect matching if there is no any pair of distinct edges in $P$ with a common vertex and any vertex in $G$ belongs to one of edges in $P$.

Let $[n] = \{1,2,\ldots,n\}$. A (finite) simplicial complex $\Delta$ on $n$ vertices, is a collection of subsets of
$[n]$ such that the following conditions hold:\\
i) $\{i\}\in\Delta$ for any $i\in [n]$,\\
ii) if $E\in\Delta$ and $F\subseteq E$, then $F\in\Delta$.\\
An element of $\Delta$ is called a face and a maximal face with respect to inclusion is called a facet. The dimension of a face $F\in\Delta$ is defined to be $|F|-1$ and dimension of $\Delta$ is maximum of dimensions of its faces. Faces with dimension 0 are called vertices.

Let $\Delta$ be a simplicial complex on
$[n]$. Let $S=K[x_1,\ldots,x_n]$ be the polynomial ring in $n$
variables with coefficients in a field $K$. Let $I_{\Delta}$
be the ideal of $S$ generated by all square-free monomials
$x_{i_1}\cdots x_{i_s}$ provided that $\{i_1,\ldots,i_s\}\not\in\Delta$. The quotient ring
$K[\Delta]=S/I_{\Delta}$ is called  Stanley-Reisner ring of the simplicial
complex $\Delta$.

Let $G$ be a graph on the vertex set $V = \{v_l,... ,v_n\}$. Let $S=K[x_1,\ldots,x_n]$. The edge ideal $I(G)$, is defined to be the ideal of $S$ generated by all square-free monomials $x_ix_j$ provided that $v_i$ is adjacent to $v_j$ in $G$. The quotient ring $R(G)=S/I(G)$ is called edge ring of $G$. We say that a set $F\subseteq V$ is an independent set in $G$ if no two of its vertices are adjacent. Define the independence complex of $G$, the simplicial complex $\Delta_G$ by
$$\Delta_G = \{F \subseteq V : F \ \mbox{is an independent set in} \ G\} .$$

Let $S$ be the polynomial ring and let $I$ be a homogeneous ideal of $S$. The depth of $S/I$, denoted by $\mbox{depth}(S/I)$, is the largest integer $r$ such that there is a sequence $f_l, \ldots, f_r$ of  homogeneous elements such that $f_i$ is not a zero-divisor in $S / (I, f_i, \ldots, f_{i-1})$ for all $1\leq i \leq r$, where $f_0$ is assumed to be $0$. Furthermore,  $(I, f_i, \ldots, f_r)\neq S$. Such a sequence is called a regular sequence. Depth is an important invariant of a ring. It is bounded by another important invariant which is Krull dimension of the ring, length of the longest chain of prime ideals. The ring $S/I$ is called Cohen-Macaulay if $\mbox{depth}(S/I) = \dim(S/I)$. A graph $G$ (a simplicial complex $\Delta$, respectively) is called Cohen-Macaulay if the ring $R(G)$ (the ring $K[\Delta]$, respectively) is Cohen-Macaulay.

A simplicial complex $\Delta$ is called pure if all its facets have the same cardinality. A graph $G$ is called unmixed if all maximal independent sets of vertices of $G$ have the same cardinality. It is clear that a graph $G$ is unmixed if and only if the simplicial complex $\Delta_G$ is pure. It is well known that a Cohen-Macaulay simplicial complex is pure, but the converse is not true, i.e., there are pure simplicial complexes which are not Cohen-Macaulay.

A pure simplicial complex $\Delta$ with vertex set $V$, is called completely balanced if there is a partition of $V$ as $C_1,\ldots,C_r$ such that each facet of $\Delta$ has exactly one vertex in common with each $C_i$. Here a partition means that $C_1\cup\cdots\cup C_r=V$ and for each $i\neq j$, $C_i\cap C_j=\varnothing$. R. Stanley has studied such simplicial complexes in \cite{St1}. He proved that, in a completely balanced simplicial complex with partition $C_1, \ldots, C_r$, the elements $\theta_1,\ldots,\theta_r$ form a homogeneous system of parameters, where
$$
\theta_i=\sum_{x\in C_i}x .
$$
Here by a homogeneous system of parameters in a standard graded ring $R$, we mean a set of homogeneous elements $\theta_1,\ldots,\theta_r$ of nonzero degrees such that $\dim(R/\langle \theta_1,\ldots,\theta_r\rangle)=0$.

R. H. Villarreal has proved in \cite{Vill2} that a bipartite graph $G$ with parts $V_1$ and $V_2$ is unmixed if and only if $|V_1|=|V_2|$ and there is an order on vertices of $V_1$ and $V_2$ as $x_1, \ldots, x_n$ and $y_1, \ldots, y_n$ respectively, such that:\\
1) $x_i\sim y_i$ for $i=1, \ldots, n$,\\
2) for each $1\leq i<j<k\leq n$ if $x_i\sim y_j$ and $x_j\sim y_k$, then $x_i\sim y_k$.

Then, M. Estrada and R. H. Villarreal in \cite{EsVil} have shown that, Cohen-Macaulayness and shellability of a bipartite graph $G$ are coincide and if $G$ is Cohen-Macaulay, then, there is a vertex $v$ in $G$ such that $G\setminus\{v\}$ is again Cohen-Macaulay.

Finally, J. Herzog and T. Hibi in \cite{HH} have proved that a bipartite graph $G$ is Cohen-Macaulay if and only if $|V_1|=|V_2|$ and there is an order on vertices of $V_1$ and $V_2$ as $x_1, \ldots, x_n$ and $y_1, \ldots, y_n$ respectively, such that:\\
1) $x_i\sim y_i$ for $i=1, \ldots, n$,\\
2) if $x_i\sim y_j$, then $i\leq j$, \\
3) for each $1\leq i<j<k\leq n$ if $x_i\sim y_j$ and $x_j\sim y_k$, then $x_i\sim y_k$.

In the above criteria, one needs to find an appropriate order on vertices of $G$ and it makes more complicated to check Cohen-Macaulayness of a given bipartite graph in practice. Here, we show that there is no need to have an order and one can check Cohen-Macaulayness of a given bipartite graph in a quite short time.

\begin{theorem}\label{main}
Let $G$ be a bipartite graph with parts $V_1$ and $V_2$. Then, $G$ is Cohen-Macaulay if and only if there is a perfect matching in $G$ as $\{x_1,y_1\}, \ldots, \{x_n,y_n\}$, such that, $x_i\in V_1$ and $y_i\in V_2$ for $i=1,\ldots,n$, and two following conditions hold.\\
1) The induced subgraph on $N(x_i)\cup N(y_i)$ is a complete bipartite graph, for $i=1,\ldots,n$. \\
2) If $x_i\sim y_j$ for $i\neq j$, then, $x_j\not\sim y_i$.
\end{theorem}

Before proving the theorem, we prove some lemmas.

\begin{lemma}\label{cm}
Let $G$ be an unmixed bipartite graph with a perfect matching $\{x_1,y_1\}, \ldots, \{x_n,y_n\}$. Then, $G$ is Cohen-Macaulay if and only if the sequence $x_1+y_1, \ldots, x_n+y_n$ is a regular sequence in $R(G)$.
\end{lemma}

\paragraph{Proof.} The sets $\{x_1,y_1\}, \ldots, \{x_n,y_n\}$ is a partition of vertices of $G$ and any maximal independent set intersects each of these sets in exactly one vertex. Thus, the simplicial complex $\Delta_{G}$ is completely balanced. By Corollary 4.2 and its Remark in \cite{St1}, $x_1+y_1, \ldots, x_n+y_n$ is a system of parameters in $R(G)$. By Theorem 17.4 in \cite{Mat} (using graded ring instead of local ring), $R(G)$ is Cohen-Macaulay if and only if every system of parameters is a regular sequence in $R(G)$.
\hfill $\Box$

\begin{lemma}\label{zdivisor}
Let $I$ be an ideal of $S=K[x_1,\ldots,x_n]$ generated by quadratic monomials. Let for some $i,j$, $1\leq i<j\leq n$, $x_i^2\not\in I$ and $x_j^2\not\in I$. Then, $\bar{x}_i+\bar{x}_j$ is zero-divisor in $S/I$ if and only if one of the following conditions hold. Here, $\bar{x}_i$ denotes the image of $x_i$ in $S/I$.\\
i) There is $x_k$, $k\not\in\{i,j\}$ such that $\bar{x}_k(\bar{x}_i+\bar{x}_j)=0$ or,\\
ii) there are integers $k,l$, $1\leq k<l\leq n$, both distinct from $i$ and $j$, such that $x_kx_l\not\in I$ and $\bar{x}_k\bar{x}_l(\bar{x}_i + \bar{x}_j)=0$.
\end{lemma}
\paragraph{Proof.}
Without loss of generality, we may assume that $i=1$ and $j=2$.
It is well known that a polynomial $f$ in $S$ belongs to a monomial ideal $I$ if and only if all monomials of $f$ are belonging to $I$. Let $\prec$ be the lexicographic order on monomials of $S$ induced by $x_1\succ x_2 \succ \cdots \succ x_n$.
Let $\bar{x}_1+\bar{x}_2$ be zero-divisor in $S/I$. Then, there is a  polynomial $h$ in $S$ such that $\bar{h}$ is nonzero in $S/I$ and $\bar{h}(\bar{x}_1+\bar{x}_2)=0$ or equivalently, $f=h(x_1+x_2)\in I$. Let $h=h_1+h_2+\cdots + h_r$ such that $h_i$'s are monomials and $h_1\succ h_2\succ\cdots\succ h_r$. We may assume that $h_1\not\in I$.  Now, $h_1x_1$ is the greatest monomial of $f$ with respect to the order $\prec$ and can not be canceled by other monomials. Therefore, $h_1x_1\in I$ and there is a quadratic monomial in generating set of $I$ which divides $h_1x_1$ and does not divide $h_1$. This monomial must be of the form $x_1x_k$ for some $k$, $1\leq k\leq n$. In other hand, $k\neq 1$ because $x_1^2\not\in I$, and $x_1\nmid h_1$ because $x_k|h_1$ and $h_1\not\in I$. According to the lexicographic order, $x_1$ does not divide any other monomial of $h$. In the polynomial $hx_1+hx_2$, in the part $hx_2$, none of monomials are divided by $x_1$. In this part, $h_1x_2$ is the greatest monomial with respect to $\prec$ and can not be canceled by other monomials and therefore, $h_1x_2\in I$. As before, there is a quadratic monomial in generating set of $I$ which divides $h_1x_2$ but not $h_1$. This monomial must be of the form $x_2x_l$ for some $2<l\leq n$. And also, $x_2\nmid h_1$. Now, $x_k|h_1$, $x_l|h_1$ and if $k=l$, then, $x_k(x_1+x_2)\in I$ and if $k\neq l$, then, $x_kx_l\not\in I$ because $x_kx_l|h_1$, and $x_kx_l(x_1+x_2)\in I$. This completes the proof in one direction. The converse is trivial.  \hfill $\Box$

\paragraph{Proof of Theorem \ref{main}.}  The proof is in 3 steps. First we prove that a bipartite graph $G$ is unmixed if and only if there is  a perfect matching in $G$ satisfying condition 1. Then, in Step 2, we prove that for an unmixed bipartite graph, condition 2 is necessary for Cohen-Macauleyness and finally in Step 3 we prove that, condition 2 is also sufficient for Cohen-Macaulayness of such a graph. \\

{\it Step 1.} Let $G$ be unmixed. There is no isolated vertex and any vertex in $V_1$ is adjacent to some vertices in $V_2$. Therefore, there is no any vertex in $V_1$ independent to the set $V_2$. This means that $V_2$ is a maximal independent set in $G$ and similarly, $V_1$ is a maximal independent set. Then, by unmixedness of $G$, $|V_1|=|V_2|$.  Let $A\subseteq V_1$ be a nonempty set and $N(A)$ be the set of all vertices in $V_2$ which are adjacent to some vertices in $A$. Suppose $|N(A)|< |A|$. There is no any edge between $A$ and $V_2\setminus N(A)$. Therefore, $A\cup (V_2\setminus N(A))$ is an  independent set and its size is strictly greater than size of $V_2$, which is a contradiction with unmixedness of $G$. Therefore, $|N(A)| \geq |A|$ for each nonempty subset $A$ of $V_1$. Therefore, by Theorem of Hall \cite{Hall}, there is a set of distinct representatives (SDR) for the set $\{\{N(v)\} : v\in V_1\}$, which determines a perfect matching between $V_1$ and $V_2$.

Now, let $V_1=\{x_1,\ldots,x_n\}$, $V_2=\{y_1,\ldots,y_n\}$ and $\{x_1,y_1\}, \ldots, \{x_n,y_n\}$ be a perfect matching in $G$. $G$ is unmixed and any maximal independent set of vertices in $G$ has cardinality $n$. Therefore, any maximal independent set intersects each edge of the perfect matching in exactly one vertex. Suppose for some $j$, $1\leq j\leq n$, the induced subgraph on $N(x_j)\cup N(y_j)$ is not complete bipartite graph. Then, there are $x\in N(y_j)$ and $y\in N(x_j)$ such that $x\not\sim y$. The set $\{x,y\}$ is independent and so there is a maximal independent set containing it. This maximal independent set does not meet the edge $\{x_j,y_j\}$ which is a contradiction. Therefore, condition 1 holds.

Conversely, let there is a perfect matching $\{x_1,y_1\}, \ldots, \{x_n,y_n\}$  in $G$ which satisfies condition 1. Let $A$ be a maximal independent set in $G$. Then $A$ meets each edge in the perfect matching in at most one vertex. Suppose that for some $j$, $1\leq j\leq n$, $A\cap \{x_j,y_j\}=\varnothing$. Then, none of $x_j$ and $y_j$ is independent to  $A$, and there are $x,y\in A$ such that $x\sim y_j$ and $y\sim x_j$. But, $x$ and $y$ are not adjacent and the induced subgraph on $N(x_j)\cup N(y_j)$ is not complete bipartite graph, which is a contradiction. Therefore, $A$ meets any edge in the perfect matching and has cardinality $n$. It means that $G$ is unmixed. \\

{\it Step 2.} Let $G$ be a bipartite graph with a perfect matching which satisfies condition 1  but condition 2 fails. That is, for some $i$ and $j$, $1\leq i<j\leq n$, we have $x_i\sim y_j$ and $x_j\sim y_i$. Then, in the quotient ring $R(G)/\langle x_i + y_i\rangle$, the element $\bar{x}_i$ is not zero and $\bar{x}_i(\bar{x}_j+\bar{y}_j)=0$ because $\bar{x}_i=-\bar{y}_i$. Therefore, $\bar{x}_j+\bar{y}_j$ is a zero-divisor in $R(G)/\langle x_i + y_i\rangle$. This means that the sequence $\bar{x}_1+\bar{y}_1, \ldots, \bar{x}_n+\bar{y}_n$ is not a regular sequence in $R(G)$ and by Lemma \ref{cm}, $R(G)$ is not Cohen-Macaulay. \\

{\it Step 3.} Let $G$ be a bipartite graph with a perfect matching satisfying condition 1. In this case, $\dim(R(G))=n$ and to prove that $R(G)$ is Cohen-Macaulay, it is enough to show that the sequence $\bar{x}_1+\bar{y}_1, \ldots, \bar{x}_n+\bar{y}_n$ is a regular sequence in $R(G)$ (Lemma \ref{cm}). For an integer $i$, $1\leq i< n$, the ring $R(G)/\langle x_1+y_1,\ldots,x_{i-1}+y_{i-1} \rangle$ can be considered to be the ring $R'(G)$ obtained by $R(G)$ with identifying variables $x_j$ with $-y_j$ for $j=1,\ldots,i-1$. By Lemma \ref{zdivisor} and its  proof, the only possibility for $\bar{x}_i+\bar{y}_i$ to be zero-divisor in $R'(G)$ is that there is $j$, $1\leq j\leq i-1$, such that $\bar{x}_j(\bar{x}_i+\bar{y}_i)=0$. Therefore, $\bar{x}_j\bar{y}_i=0$ and $\bar{x}_j\bar{x}_i=0$ or equivalently, $\bar{y}_j\bar{x}_i=0$. Therefore, $x_j\sim y_i$ and $y_j\sim x_i$. But, in this case, condition 2 fails. This completes the proof. \hfill $\Box$

\begin{proposition}
Condition 1 in Theorem 1 which is equal to unmixedness of a bipartite graph  is also equal to saying that non of the polynomials $x_1+y_1, \ldots, x_n+y_n$ are zero-divisor in $R(G)$.
\end{proposition}

\paragraph{Proof.}  It is clear by Lemma \ref{zdivisor} and Theorem \ref{main}. \hfill $\Box$

\begin{remark}\label{connect} Condition 2 in Theorem 1 is equal to say that, for each $i$ and $j$, $1\leq i<j\leq n$, the induced subgraph on vertices $\{x_i,y_i,x_j,y_j\}$ has connected complement.
\end{remark}

\begin{corollary}
Let $G$ be a bipartite Cohen-Macaulay graph and $\{x_i,y_i\}$ be any edge in the perfect matching mentioned in Theorem 1. Then, $G\setminus\{x_i,y_i\}$ is again Cohen-Macaulay.
\end{corollary}
\paragraph{Proof.} Here, by $G\setminus\{x_i,y_i\}$ we mean the graph obtained by deleting vertices $x_i$ and $y_i$ and all edges passing through one of these vertices. It is clear that if condition 1 or 2 in Theorem 1 holds for $G$, then, it holds for $G\setminus\{x_i,y_i\}$ for each $i=1,\ldots,n$.
\hfill $\Box$

\begin{proposition}
Let $G$ be a bipartite Cohen-Macaulay graph with parts $V_1$ and $V_2$. Then, there is at least one vertex of degree one in each part.
\end{proposition}
\paragraph{Proof.} Let $y$ be a vertex in $V_2$ such that for any other vertex $y'\in V_2$, we have $\deg(y')\leq\deg(y)$. Let $x\in V_1$ be the vertex such that $\{x,y\}$ is in a perfect matching in $G$. We have $\deg(x)\geq 1$. If $\deg(x)>1$, then there is a vertex $y'\in V_2\setminus \{y\}$ such that $x\sim y'$. Let $x'$ be a vertex in $V_1\setminus\{x\}$ such that $\{x',y'\}$ is in the perfect matching. $G$ is Cohen-Macaulay then, the induced subgraph on  $N(x)\cup N(y)$ is a complete bipartite graph and $x'\not\in N(y)$. Then, $y'$ is adjacent to each vertex in $N(y)\cup \{x'\}$. Therefore, $\deg(y')>\deg(y)$ which is a contradiction. Therefore, $\deg(x)=1$. \hfill $\Box$

Let $G$ be a Cohen-Macaulay bipartite graph. There are some vertices in both parts with degree one. If we remove the vertex adjacent to a vertex of degree one, the edge consisting these two vertices in a perfect matching will be removed and the remaining graph is also Cohen-Macaulay. 

\begin{corollary}\label{unique}
Let $G$ be a Cohen-Macaulay bipartite graph. There is a unique perfect matching in $G$.
\end{corollary}
\paragraph{Proof.} By Theorem 1, there is a perfect matching. Let $V_1$ and $V_2$ be two parts of $G$. Let $P$ be a perfect matching in $G$. By the above proposition, there is a vertex of degree one in $V_1$. Let $x_1$ be the vertex and $y_1\in V_2$ be the unique vertex adjacent to $x_1$. Then $\{x_1,y_1\}\in P$. The graph $G\setminus\{x_1,y_1\}$ is again Cohen-Macaulay and $V_1\setminus\{x_1\}$ has a vertex of degree one as $x_2$. Let $y_2\in V_2\setminus\{y_1\}$ be the unique vertex adjacent to $x_2$. Then, $\{x_2,y_2\}\in P$. Continuing this process, determines $P$ uniquely.
\hfill $\Box$

\begin{corollary}
Let $G$ be an unmixed bipartite graph. Then, the following conditions are equivalent.\\
i) $G$ is Cohen-Macaulay. \\
ii) There is a unique perfect matching in $G$.\\
iii) For each two edges $e_1, e_2$ in a perfect matching, complement of the induced subgraph on vertices of $e_1$ and $e_2$ is connected.
\end{corollary}
\paragraph{Proof.} (i$\to$ii) is proved in Corollary \ref{unique}. Let $G$ be unmixed but not Cohen-Macaulay. Then, there is a perfect matching and two edges in the perfect matching as $\{x_i,y_i\}$ and $\{x_j,y_j\}$ such that $x_i\sim y_j$ and $x_j\sim y_i$. Substituting $\{x_i,y_i\}$ and $\{x_j,y_j\}$ by $\{x_i,y_j\}$ and $\{x_j,y_i\}$, we get a different perfect matching. This proves (ii$\to$i). Equality of i and iii is clear by Theorem \ref{main} and Remark \ref{connect}. \hfill $\Box$

For a given bipartite graph $G$, we present a fast polynomial-time algorithm to check wether $G$ is Cohen-Macaulay or not.

\begin{algorithm}
Let $G$ be a given bipartite graph with $m$ vertices.
\begin{itemize}
\item[]{\em Step~1.} Take $i=0$. If $m$ is not even, then, go to Step~7.
\item[]{\em Step~2.} If there is no any vertex with degree 1 in $G$, go to Step~7.
\item[]{\em Step~3.} Take $i=i+1$. Choose a vertex of degree one and name it $x_i$. Name the vertex adjacent to $x_i$ to be $y_i$. Take $G=G\setminus\{x_i,y_i\}$. If $i<n$, go to Step~2.
\item[]{\em Step~4.} If there is $j$, $1\leq j\leq n$ such that, a vertex in $N(x_j)$ and a vertex in $N(y_j)$ are not adjacent, then,  go to Step~7.
\item[]{\em Step~5.} If there are $i, j$, $1\leq i<j\leq n$ such that $x_i\sim y_j$ and $x_j\sim y_i$, then, go to Step~7.
\item[]{\em Step~6.} Write "G is Cohen-Macaulay" and end the algorithm.
\item[]{\em Step~7.} Write "G is not Cohen-Macaulay" and end the algorithm.
\end{itemize}
\end{algorithm}

In Step 3, $G\setminus\{x_i,y_i\}$ is the induced subgraph of $G$ on vertex set $V(G)\setminus\{x_i,y_i\}$.

Note that the assumption that there is no  vertex of degree zero in $G$ is not really a restriction in the class of all bipartite graphs for Cohen-Macaulayness. Because, any graph with only one vertex is Cohen-Macaulay and disjoint union of two graphs is Cohen-Macaulay if and only if both of them are Cohen-Macaulay. Therefore, in a given bipartite graph $G$ we may omit all isolated vertices and check Cohen-Macaulayness of the remaining graph.

Some of results of this paper were already known. For example, equality of unmixedness of a bipartite graph with  condition 1 in Theorem \ref{main} is proved in \cite{Riv}. But, the aim of this work was gathering together these results and reformulate and reprove them in a constructive way such that an algorithm can be obtained. Also we hope that the proofs in this paper give some ideas to find the same results in a larger class consisting of $r$-partite graphs which have some separated maximal cliques covering all vertices.

\end{document}